\documentclass[11pt]{article}
\usepackage{amsmath,amsfonts,amssymb,amsthm,fancyhdr,bm,mathtools,enumitem}
\usepackage[hidelinks]{hyperref}
\usepackage{fullpage}
\usepackage[dvipsnames]{xcolor}
\usepackage[capitalize]{cleveref}
\usepackage[color=Purple!50!white]{todonotes}
\usepackage[numbers,sort]{natbib}
\usepackage{tikz,ifthen}

\usetikzlibrary{decorations.markings,arrows.meta}
\tikzset{vert/.style={draw, fill=black, circle, inner sep=2pt}}

\newtheorem{theorem}{Theorem}%[section]
\newtheorem{lemma}[theorem]{Lemma}

\newtheorem{qu}[theorem]{Open problem}
\theoremstyle{definition}

\crefname{equation}{equation}{equations}
\crefname{lemma}{Lemma}{Lemmas}
\crefname{proposition}{Proposition}{Propositions}
\crefname{claim}{Claim}{Claims}
\crefname{theorem}{Theorem}{Theorems}
\crefname{conjecture}{Conjecture}{Conjectures}
\crefname{figure}{Figure}{Figures}

\newlist{lemenum}{enumerate}{1}
\setlist[lemenum]{label=(\alph*), ref=\thelemma(\alph*)}
\crefalias{lemenumi}{lemma}

\let\leq\leqslant

\let\geq\geqslant

\title{Infinitely many minimally non-Ramsey size-linear graphs}
\author{Yuval Wigderson\thanks{Institute for Theoretical Studies, ETH Z\"urich, 8006 Z\"urich, Switzerland. 
Supported by Dr.\ Max R\"{o}ssler, the Walter Haefner Foundation, and the ETH Z\"{u}rich Foundation. Email: {\tt{yuval.wigderson@eth-its.ethz.ch}}}}
\date{}

\begin{document}
\maketitle
\begin{abstract}
	A graph $G$ is said to be Ramsey size-linear if $r(G,H) =O_G (e(H))$ for every graph $H$ with no isolated vertices. Erd\H os, Faudree, Rousseau, and Schelp observed that $K_4$ is not Ramsey size-linear, but each of its proper subgraphs is, and they asked whether there exist infinitely many such graphs. In this short note, we answer this question in the affirmative.
\end{abstract}
Given two graphs $G,H$, their \emph{Ramsey number} $r(G,H)$ is the least integer $N$ such that every two-coloring of $E(K_N)$ contains a monochromatic copy of $G$ in the first color, or of $H$ in the second color. Our understanding of $r(G,H)$ is rather limited in general, but a great deal is known in certain special cases. For example, Chv\'atal \cite{MR465920} proved that\footnote{We use $v(H)$ and $e(H)$ to denote the number of vertices and edges, respectively, of a graph $H$.} $r(T,K_n)=(v(T)-1)(n-1)+1$ for every tree $T$, and Sidorenko \cite{MR1223692} proved that $r(K_3,H) \leq 2e(H)+1$ for every graph $H$ with no isolated vertices, which is tight if $H$ is a tree or a matching.

Generalizing this second example, Erd\H os, Faudree, Rousseau, and Schelp \cite{MR1264714} defined a \emph{Ramsey size-linear graph} to be a graph $G$ for which $r(G,H) \leq C_G\cdot e(H)$ for every graph $H$ with no isolated vertices, where $C_G>0$ is a constant depending only on $G$. Thus, Sidorenko's result \cite{MR1223692} implies that $K_3$ is Ramsey size-linear. On the other hand, $K_4$ is not Ramsey size-linear, since $r(K_4,K_n)=\omega(n^2)$ \cite{MR491337,MR4713025}, whereas $K_n$ has $\binom n2 = O(n^2)$ edges.

Erd\H os, Faudree, Rousseau, and Schelp \cite{MR1264714} observed that in fact, $K_4$ is minimally non-Ramsey size-linear, in the sense that every proper subgraph of $K_4$ is Ramsey size-linear. They asked whether there exist infinitely many such graphs, or even more restrictively, whether there exist any examples besides $K_4$. This question was reiterated in \cite{MR1370501,MR1425200}, and appears as problem 79 on Bloom's Erd\H os problems website \cite{erdosproblems}. In this note, we show that there are infinitely many such graphs.
\begin{theorem}\label{thm:main}
	There exist infinitely many graphs $G$ which are not Ramsey size-linear, but every proper subgraph $G' \subsetneq G$ is Ramsey size-linear.
\end{theorem}
In the course of the proof of \cref{thm:main}, we shall need the following three simple facts.
\begin{lemma}\label{lem:forests}
	Every forest is Ramsey size-linear.
\end{lemma}
Indeed, it suffices to prove this for trees, since every forest is a subgraph of a tree. Every graph $H$ with no isolated vertices is a subgraph of $K_{2e(H)}$, so \cref{lem:forests} follows immediately from the result of Chv\'atal \cite{MR465920} mentioned above. Substantially stronger results than \cref{lem:forests} are proved in \cite[Theorems 3--5]{MR1264714}.

\begin{lemma}[{\cite[Corollary 1]{MR1264714}}]\label{lem:avg deg}
	If $e(G) \geq 2v(G)-2$, then $G$ is not Ramsey size-linear.
\end{lemma}
Indeed, using the Lov\'asz local lemma, one can show that $r(G,K_n) = \Omega((n/{\log n})^{\frac{e(G)-1}{v(G)-2}})$ (see \cite{MR491337,MR1264714} for details). If $e(G) \geq 2v(G)-2$ then this exponent is strictly greater than $2$, hence $K_n$ witnesses that $G$ is not Ramsey size-linear.
\begin{lemma}\label{lem:girth}
	For every $g\geq 3$, there exists a graph with girth at least $g$ and average degree at least $4$.
\end{lemma}
The existence of such a graph follows immediately from a standard probabilistic deletion argument. For explicit constructions, one can use the Ramanujan graphs of Lubotzky--Phillips--Sarnak \cite{MR963118}, for example. One can also greedily construct such graphs by repeatedly joining pairs of vertices at distance at least $g$, while controlling the maximum degree; see e.g.\ \cite[Theorem 2.13]{2403.13571} for details.

With these preliminaries, we are ready to prove \cref{thm:main}.
\begin{proof}[Proof of \cref{thm:main}]
	Suppose for contradiction that there exist only finitely many such graphs, say $G_1,\dots,G_k$. By \cref{lem:forests}, each $G_i$ contains at least one cycle, say of length $\ell_i$. Let $g = 1+\max \{\ell_1,\dots,\ell_k\}$. By \cref{lem:girth}, there exists a graph $G_0$ with girth at least $g$ and average degree at least $4$.
	Note that no $G_i$ is a subgraph of $G_0$, since $G_i$ has a cycle of length $\ell_i$ but $G_0$ does not.

	Moreover, since the average degree of $G_0$ is at least $4$, we have $e(G_0)\geq 2v(G_0)$, hence $G_0$ is not Ramsey size-linear by \cref{lem:avg deg}. Let $G$ be an inclusion-wise minimal subgraph
	of $G_0$ which is not Ramsey-size linear. By construction, $G$ is not Ramsey size-linear, but every proper subgraph of it is. Moreover, $G \notin \{G_1,\dots,G_k\}$, since $G$ is a subgraph of $G_0$ but none of $G_1,\dots,G_k$ is. This contradiction completes the proof.
\end{proof}

We remark that this proof is non-constructive, in the sense that it does not supply any example of a minimally non-Ramsey size-linear graph. As such, the following natural problem remains open.
\begin{qu}
	Give an example of a minimally non-Ramsey size-linear graph other than $K_4$.
\end{qu}
The proof of \cref{thm:main} implies that if one starts with a $K_4$-free graph with average degree at least $4$, such as $K_{2,2,2}$ or $K_{4,4}$, then some subgraph of it is minimally non-Ramsey size-linear, but it seems difficult to identify such a subgraph.

\vspace{2pt}

\noindent\textbf{Acknowledgments:} I am grateful to Domagoj Brada\v c and Jacob Fox for discussions and comments on an earlier draft, and to the anonymous referees for their helpful suggestions.

% \bibliographystyle{yuval}
% \bibliography{refs.bib}

\end{document}